\documentclass[pdflatex,sn-mathphys-num]{sn-jnl}% Math and Physical Sciences Numbered Reference Style

\usepackage[T1]{fontenc}
\usepackage[utf8]{inputenc}
\usepackage{graphicx}%
\usepackage{multirow}%
\usepackage{amsmath,amssymb,amsfonts}%
\usepackage{amsthm}%
\usepackage{mathrsfs}%
\usepackage[title]{appendix}%
\usepackage{xcolor}%
\usepackage{textcomp}%
\usepackage{manyfoot}%
\usepackage{booktabs}%
\usepackage{algorithm}%
\usepackage{algorithmicx}%
\usepackage{algpseudocode}%
\usepackage{listings}%
\allowdisplaybreaks
\theoremstyle{thmstyleone}%
\newtheorem{theorem}{Theorem}% 
\newtheorem{proposition}[theorem]{Proposition}% 
\newtheorem{lemma}[theorem]{Lemma}

\theoremstyle{thmstyletwo}%
\newtheorem{example}{Example}%

\theoremstyle{thmstylethree}%

\renewcommand{\epsilon}{\varepsilon}

\begin{document}

\title[Torus Bifurcation in the Hopf-Langford system through Averaging Theory]{Torus Bifurcation in the Hopf-Langford type system through Averaging Theory}

\author*[1]{\fnm{Gustavo} \sur{Domingues}}\email{marra@unifei.edu.br}

\affil[1]{\orgdiv{Instituto de Ci\^encias Puras e Aplicadas}, \orgname{Universidade Federal de Itajub\'a, \textit{campus} Itabira}, \orgaddress{\street{Rua Irm\~a Ivone Drummond, 200, Distrito Industrial II}, \city{Itabira}, \postcode{35903-087}, \state{Minas Gerais}, \country{Brasil}}}

\abstract{We study the existence of periodic solutions and invariant tori bifurcating from a zero-Hopf equilibrium in a Hopf-Langford type system. Using second-order averaging theory, we establish explicit parameter conditions under which the system admits a periodic solution near the origin, together with its stability. We then apply recent results relating invariant tori to Neimark-Sacker bifurcations of the associated Poincar\'e map to obtain a smooth bifurcation curve along which a unique invariant torus surrounds the periodic solution. The relevant first Lyapunov coefficient is computed explicitly as a power series in the perturbation parameter, and we show that, under the conditions considered, this coefficient is always positive; consequently, whenever it exists, the bifurcating torus is unstable, regardless of the specific parameter values chosen within our hypotheses. This provides a structural obstruction to the existence of asymptotically stable tori in this family of systems. We illustrate our results with two numerical examples, one exhibiting an isolated asymptotically stable periodic solution with no bifurcating torus, and another exhibiting a periodic solution surrounded by an unstable invariant torus, the latter also illustrated through its Poincar\'e map.}



\maketitle

\section{Introduction}\label{sec1}

In 1948 E. Hopf \cite{hopf1948mathematical} studied turbulence phenomena from a dynamical point of view, using a system of equations with structural characteristics analogous to the Navier-Stokes equations. Such a system allows an infinite sequence of bifurcations, each one adding an independent period to a quasi-periodic motion. It was used by Langford in 1979 \cite{langford1979periodic} to study bifurcation of equilibria and periodic solutions and is given as below:
\begin{equation}\label{langfordoriginal}
\left\{
\begin{array}{rcl}
\dot{x}&=&(\nu -1 )x-y +x z\\
\dot{y}&=&x+(\nu -1 )y+y z\\
\dot{z}&=&\nu  z- \left(x^2+y^2+z^2\right)
\end{array}
\right.,
\end{equation}
where $\nu$ is a real parameter. Langford conjectured the existence of invariant tori in this system. This system was later studied in \cite{hassard1981theory} where the authors prove successive local Hopf bifurcations and provide a global geometric description for the emergence of an invariant torus.

  Several differential systems inspired by \eqref{langfordoriginal}, called Generalized Hopf-Langford (GHL) systems, are studied accounting for different aspects of their dynamics, including the existence of invariant tori (see \cite{langford1984numerical, su2008bifurcation,guo2017steady,yang2018complex,
bashkirtseva2020stochastic,nikolov2021assessing,nikolov2021completely,
fu2021bifurcations,musafirov2022admissible,
zhong2023invariant,guo2023stability,yumagulov2023langford,vassilev2024first}).

 A specific system belonging to this class is the one given below (which we call Hopf-Langford type as there are no added terms):
\begin{equation}\label{langfordsystem}
\left\{
\begin{array}{rcl}
\dot{x}&=&x (\mu -\alpha )-\beta  y +x z\\
\dot{y}&=&\beta  x+y (\mu -\alpha )+y z\\
\dot{z}&=&\mu  z-\gamma \left(x^2+y^2+z^2\right)
\end{array}
\right.,
\end{equation}
with all parameters real and the derivatives of $x,y,z$ are taken with respect to time. 

Our approach in this work studies an invariant torus bifurcating from a \textit{zero-Hopf equilibrium}: in a three-dimensional autonomous differential system it is an equilibrium point whose eigenvalues are 0 and $\pm i \, \omega $ ($\omega>0$). Zero-Hopf equilibria are well-known indicators of rich local dynamical behavior, and we remark that there is no general theory that completely describes the dynamics around a zero-Hopf equilibrium.

The classical \emph{averaging method} is one of the main tools used to detect periodic solutions in non-autonomous periodic differential equations and systems. It consists of generating a sequence of functions $\textbf{g}_i$ (called the \textit{i}-th \textit{order averaged function}) for which their simple zeros correspond to isolated periodic solutions for the differential system.

Concerning invariant tori in differential systems, in works \cite{elsamaloty1984averaging, hale1961integral, marsden2012hopf, sacker1965new} there are several correlations between the existence of invariant tori and zero-Hopf bifurcations on an averaged system for a certain class of differential equations. In \cite{candido2020torus}, authors associate the bifurcation of an invariant torus on a non-autonomous differential systems with a Neimark-Sacker bifurcation of the corresponding Poincar\'e map, using averaging theory. Generic conditions over the averaged functions are established to guarantee the existence of a torus bifurcation curve in a system with two parameters. This analysis is also extended to higher-order averaged function. Results from \cite{candido2020torus} have been used to study periodic solutions and invariant tori bifurcating from zero-Hopf equilibria in several distinct differential systems \cite{candido2020periodic, candido2022zero, candido2024invariant, candido2024analytical}.

In this work, we study the existence of periodic orbits and invariant torus bifurcating from a zero-Hopf equilibrium on system \eqref{langfordsystem}. We use recent results from averaging theory to provide generic conditions ensuring the existence of an invariant torus bifurcating from a periodic orbit via a secondary Hopf (Neimark-Sacker) bifurcation. In particular, we establish parameter conditions that ensures an invariant torus bifurcating from the origin of system \eqref{langfordsystem}. 

If $\alpha = 0$, $\mu = 0$ and $\beta=\omega > 0$, it is straightforward to check that under such conditions system \eqref{langfordsystem} has a zero-Hopf equilibrium at the origin. Consider the degree 2 regular perturbation of parameters given by
\begin{equation}\label{perturbation}
\begin{array}{rcl}
\alpha &=& \alpha _1 \epsilon + \alpha_2 \epsilon^2+ \mathcal{O}(\epsilon^{3}),\\
\beta &=& \omega+\beta _1 \epsilon + \beta_2 \epsilon^2 + \mathcal{O}(\epsilon^{3}),\\
\gamma &=& \gamma _0+\gamma _1 \epsilon+  \gamma_2 \epsilon^2 +\mathcal{O}(\epsilon^{3}),\\
\mu  &=&\mu _1 \epsilon +  \mu_2 \epsilon^2+ \mathcal{O}(\epsilon^{3}).
\end{array}
\end{equation}

Consider the following set of conditions:
\begin{equation}\label{conditions} 
\begin{array}{l}
 \quad \alpha_1 \neq 0,\\ \\
  \quad \gamma_0>0,\\ \\
  (\mu_1 - \alpha_1) (\alpha_1 \gamma_0 - (\gamma_0 +1)\mu_1) > 0,\\ \\
  4 \alpha _1^2 \gamma _0 \left(\gamma _0+2\right)-4 \alpha _1 \left(\gamma _0+2\right) \left(2 \gamma _0+1\right) \mu _1+\left(2 \gamma _0+3\right){}^2 \mu _1^2 <0.
  \end{array}
\end{equation}

The following theorem is similar to a result from \cite{yang2018complex} in a slightly different setting.

\begin{theorem}\label{teoremaprincipal0}
Let $(\alpha, \beta, \gamma, \mu)$ be as in \eqref{perturbation} on system \eqref{langfordsystem}. Assume that conditions \eqref{conditions} hold. Then, for $\epsilon$ sufficiently small, system \eqref{langfordsystem} admits a periodic solution $\phi(t,\epsilon)$ satisfying $\phi(t,\epsilon) \rightarrow 0$ as $\epsilon \rightarrow 0$. Moreover, this periodic solution is asymptotically stable (resp. unstable) if
$$
2 \gamma _0 (\mu_1 - \alpha_1)+\mu _1 <0 \text{ (resp. } >0 ).
$$

\end{theorem}

Let
\begin{equation}\label{omega0}
\omega_0 = \dfrac{\alpha _1 \sqrt{ 2 \gamma _0}}{\omega(2 \gamma _0 +1) },
\end{equation}
\begin{equation}\label{mu0}
\widehat{\mu}_0 = \dfrac{2 \alpha _1 \gamma _0}{2 \gamma _0+1},\end{equation}
\begin{equation}\label{mu1}
\begin{array}{rcl}
\widehat{\mu}_1 &=& \dfrac{1}{\omega\left(2 \gamma _0+1\right){}^3 } \left( -4 \gamma _0^2 \omega  \left(-4 \pi  \alpha _2+6 \pi  \mu _2+\omega  \omega_0^2\right)+\gamma _0 \left(8 \pi  \alpha _1 \gamma _1 \omega \right. \right.\\ \\
&&  \left. -4 \omega  \left(-\pi  \alpha _2+3 \pi  \mu _2+\omega  \omega_0^2\right)+8 \pi ^2 \alpha _1^2\right)-\omega  \left(-4 \pi  \alpha _1 \gamma _1+2 \pi  \mu _2+\omega  \omega_0^2\right) \\ \\
&& \left.+16 \pi  \gamma _0^3 \omega  \left(\alpha _2-\mu _2\right) \right),\\ \\
\end{array}
\end{equation}
\begin{equation}\label{mu2}
\begin{array}{rcl}
\widehat{\mu}_2 &=& \dfrac{1}{\left(2 \gamma _0+1\right)^5 \omega ^2} \left( 16 \pi ^3 \alpha _1^2 \gamma _0 \left(2 \gamma _0+1\right) \left(\beta _1 \left(2 \gamma _0+1\right)-6 \gamma _1 \omega \right)+64 \pi ^4 \alpha _1^3 \gamma _0  \right.\\ \\
&& -8 \alpha _1 \left(2 \pi  \gamma _0+\pi \right){}^2 \omega  \left(-4 \pi  \alpha _2 \gamma _0-\beta _1 \left(2 \gamma _0+1\right) \gamma _1+\omega  \left(2 \gamma _1^2+\omega_0^2\right)\right) \\ \\
&& -2 \left(2 \gamma _0+1\right){}^3 \omega+\pi  \beta _1 \left(2 \gamma _0+1\right) \left(2 \pi  \left(2 \gamma _0 \left(\mu _2-\alpha _2\right)+\mu _2\right)+\omega  \omega_0^2\right) \\ \\
&&   +\omega  \left(\left(2 \gamma _0+1\right) \omega  \omega_0 \right. \left. -2 \pi  \gamma _1 \left(2 \pi  \alpha _2+\omega  \omega_0^2\right)\right),\\ \\
\end{array}
\end{equation}
\begin{equation}\label{l12}
\ell_{1,2} = \dfrac{\left(4 \pi ^2-2\right) \gamma _0^2+3 \left(4 \pi ^2-1\right) \gamma _0+\pi ^2}{16 \omega ^2}.
\end{equation}

Note that, under the assumption $\gamma_0>0$, we have $\ell_{1,2}>0$. We state our main theorem below.

\begin{theorem}\label{teoremaprincipal}
Let $(\alpha, \beta, \gamma, \mu)$ be as in \eqref{perturbation} on system \eqref{langfordsystem}. Assume that conditions \eqref{conditions} hold.
\begin{itemize}
\item[(a)]  For a sufficiently small $\epsilon$, the system has a periodic solution $\phi(t,\epsilon)$ satisfying $\phi(t,\epsilon) \rightarrow 0$ when $\epsilon \rightarrow 0$.

\item[(b)] Let $\phi(t,\epsilon) = \phi(t,\mu_1,\epsilon)$ be the periodic solution from item \textit{(a)}. Then there exists a smooth curve $\mu(\epsilon) = \widehat{\mu_0} + \widehat{\mu}_1 \epsilon + \widehat{\mu_2} \epsilon^2+\mathcal{O}(\epsilon^3)$ and intervals $J_\epsilon$ containing $\mu(\epsilon)$ such that a unique invariant torus bifurcates from the periodic solution $\phi(t,\mu(\epsilon),\epsilon)$ as $\mu_1$ passes through $\mu(\epsilon)$. The torus exists for $\mu \in J_\epsilon$, provided $\ell_{1,2} (\mu_1- \mu(\epsilon)) <0$, and surrounds the periodic solution. Moreover, since $\ell_{1,2}>0$ the torus is unstable whenever it exists.
\end{itemize}
\end{theorem}

The paper is organized as follows. In Section 2 we present the fundamentals of the averaging methods and torus bifurcation in this setting. In Section 3 we prove our main theorem. Section 4 is dedicated to numerical examples.

\section{Averaging Theory for periodic solutions and invariant tori}

	We divide this section into two parts: one concerning the averaging method and another devoted to the tools needed to detect bifurcating tori. We refer the reader to \cite{sanders2007averaging} and \cite{candido2020torus} and references therein for a deeper coverage of these subjects.

\subsection{Averaging Theory and Existence of periodic solutions}

We limit our presentation to second order averaging. Consider a non-autonomous perturbed differential system in the following standard form:
\begin{equation}\label{eqavrnormalform}
\dot{\mathbf{x}} =  \epsilon \mathbf{F}_1(t,\mathbf{x}) + \epsilon^2 \mathbf{F}_2(t,\mathbf{x})+\epsilon^{3}\widetilde{\mathbf{F}_3}(t,\mathbf{x},\epsilon),
\end{equation}
where $t\in\mathbb{R}$, $\mathbf{x}\in \Omega$, an open subset of $\mathbb{R}^n$, and $| \epsilon | < \epsilon_0$ for some sufficiently small $\epsilon_0>0$. Assume that the functions $\mathbf{F}_1$, $\mathbf{F}_2$ and $\widetilde{\mathbf{F}_3}$ are sufficiently smooth and $T$-periodic in $t$. This periodicity allows system \eqref{eqavrnormalform} to be viewed as a vector field on the cylinder
\[
(t,\mathbf{x}) \in \mathbb{S}^1 \times \Omega,
\quad \mathbb{S}^1 = \mathbb{R}/T\mathbb{Z}.
\]

Let $\mathbf{x}(t,\mathbf{z},\epsilon)$ be the solution of \eqref{eqavrnormalform} that satisfies $\mathbf{x}(0,\mathbf{z},\epsilon)=\mathbf{z}$. Expressing it as a power series in $\epsilon$ around $\epsilon=0$ one can write it as
$$
\mathbf{x}(t,\mathbf{z},\epsilon) = \mathbf{z} + \epsilon \mathbf{y}_1(t,\mathbf{z}) + \epsilon^2 \dfrac{\mathbf{y}_2 (t,\mathbf{z})}{2} + \mathcal{O}(\epsilon^3).
$$
Using this, the Poincar\'e map $\Pi(\mathbf{z},\epsilon )= \mathbf{x}(T,\mathbf{z},\epsilon)$ can be expressed using series as
$$
\Pi(\mathbf{z},\epsilon) = \mathbf{z} + \epsilon \mathbf{g}_1 (\mathbf{z}) + \epsilon^2 \mathbf{g}_2(\mathbf{z}) + \mathcal{O}(\epsilon^3).
$$
The functions $\mathbf{g}_i$ are called the \emph{averaged functions of order $i$} of system \eqref{eqavrnormalform}. Expressions for averaged functions of higher order can be found in \cite{candido2017persistence} and \cite{novaes2017equivalent}. In particular, we have: 

\begin{equation}\label{g1geral}
\mathbf{g}_1(\mathbf{z}) = \displaystyle\int_0^T \mathbf{F}_1(\tau,\mathbf{z}) d\tau.
\end{equation}
and
\begin{equation}\label{g2geral}
\mathbf{g}_2(\mathbf{z}) = \displaystyle\int_0^T \mathbf{F}_2(\tau,\mathbf{z}) +  \dfrac{\partial \mathbf{F}_1}{\partial \mathbf{z}} (\tau, \mathbf{z}) \mathbf{y}_1(\tau, \mathbf{z})  d\tau.
\end{equation}

In this setting, we have:
\begin{proposition}[\cite{llibre2014erratum}]\label{zerodapromed} Suppose that there exists $\mathbf{x}_0 \in \Omega$ such that $\mathbf{g}_1(\mathbf{x}_0)=0$ and $\det D\mathbf{g}_1 (\mathbf{x}_0) \neq 0$. Then, for $\epsilon$ sufficiently small the system \eqref{eqavrnormalform} admits a unique isolated $T$-periodic solution $\phi(t,\epsilon)$ satisfying
$\phi(0,\epsilon) \to \mathbf{x}_0$ as $\epsilon \to 0$. Moreover, the stability of this periodic solution is determined by the eigenvalues of the Jacobian matrix $D\mathbf{g}_1(\mathbf{x}_0)$: if all eigenvalues have negative (resp. positive) real parts, the periodic orbit is asymptotically stable (resp. unstable). 
\end{proposition}

We remark that, up to a constant factor of $\frac{1}{T}$, the definitions of $\mathbf{g}_1$ and $\mathbf{g}_2$ coincide with the usual first order and second order averaged functions.

\subsection{Bifurcation of Invariant tori}

In this section we outline the method presented in \cite{candido2020torus}, in which generic conditions are established over the averaged functions of a system \eqref{eqavrnormalform2} to detect a bifurcation of an invariant closed curve around a fixed point on a transversal section $\Sigma$ as the fixed point changes stability (a \emph{Neimark-Sacker} bifurcation, in discrete dynamical systems). This invariant curve in $\Sigma$ corresponds to an invariant torus on the original system. The Neimark-Sacker bifurcation has a well-known normal form which is described in \cite{kuznetsov1998elements}.

Consider system \eqref{eqavrnormalform} with an additional parameter $\mu$:

\begin{equation}\label{eqavrnormalform2}
\dot{\mathbf{x}} =  \epsilon \mathbf{F}_1(t,\mathbf{x},\mu)+\epsilon^2 \mathbf{F}_2(t,\mathbf{x},\mu)+\epsilon^{3}\widetilde{\mathbf{F}_3}(t,\mathbf{x},\epsilon,\mu),
\end{equation}
with $(t,\mathbf{x},\epsilon,\mu) \in \mathbb{R} \times \Omega \times (-\epsilon_0,\epsilon_0) \times J$, $\Omega$ an open subset of $\mathbb{R}^2$ and $J$ an open interval.

	Let the first and second order averaged functions of \eqref{eqavrnormalform2} be $\mathbf{g}_1(\mathbf{x},\mu) $ and $\mathbf{g}_2(\mathbf{x},\mu)$ where ${\bf x} = (x,y)$. The Poincar\'e map induced by this particular solution may be written as
	\begin{equation}\label{poincaremap}
	\mathbf{x}\mapsto P(\mathbf{x},\epsilon,\mu) = \mathbf{x}(T,\mathbf{z},\epsilon,\mu) = \mathbf{x} + \epsilon \mathbf{G}( \mathbf{x}, \mu, \epsilon)
	\end{equation}
where
\begin{equation}\label{poincaremapG}
\mathbf{G}(\mathbf{x},\mu,\epsilon) = \mathbf{g}_1(\mathbf{x},\mu) + \epsilon \mathbf{g}_2 (\mathbf{x},\mu) + \epsilon^2 \widetilde{\mathbf{G}}(\mathbf{x},\mu,\epsilon)
\end{equation}
on the transversal section $\Sigma = \{t=0\}\times \Omega$.  

 Consider the \textit{first order averaged system} of \eqref{eqavrnormalform2}, given by
\begin{equation}\label{sistemaaveraged}
\dot{\bf{x}} = \epsilon \, \bf{g}_1(\bf{x},\mu).
\end{equation}

Assume that system \eqref{sistemaaveraged} has a Hopf point $(\bf{x}_{\mu_0},\mu_0)$: that is, it verifies $\bf{g}_1(\bf{x}_{\mu_0},\mu_0) = 0$ and the Jacobian matrix $D \bf{g}_1(\bf{x}_{\mu_0},\mu_0)$ has eigenvalues $\pm i \omega_0$ for some positive $\omega_0$. By using the Implicit Function Theorem, one obtains a continuous curve $\mu \to \bf{x}_\mu \in \Omega$, $\mu \in J$ such that $\bf{g}_1(\bf{x_\mu},\mu)=0$. The pair of complex conjugate eigenvalues $a(\mu) \pm i b(\mu)$ of $D \bf{g}_1 (\bf{x}_\mu,\mu)$ verifies $a(\mu_0)=0$ and $b(\mu_0)=\omega_0$. Therefore, the existence of a Hopf point $(\bf{x}_{\mu_0},\mu_0)$ for system \eqref{sistemaaveraged} is equivalent to the following hypothesis:

\begin{itemize}
\item[\bf{H1.}] \emph{There exists a continuous curve $\mu \to \bf{x}_\mu$, $\mu(0)=\mu_0 \in J$ such that $\bf{g}_1(\bf{x_\mu},\mu)=0$ for every $\mu \in J$. The pair of complex conjugate eigenvalues $a(\mu) \pm i b(\mu)$ of $D \bf{g}_1 (\bf{x}_\mu,\mu)$ verifies $a(\mu_0)=0$ and $b(\mu_0)=\omega_0$. Without loss of generality (under an adequate change of coordinates), we assume $D\bf{g}_1(\bf{x}_{\mu_0},\mu_0)$ is in its real Jordan normal form.}
\end{itemize}

One can find open intervals $J_0$ containing $\mu_0$ and $0<\epsilon_1<\epsilon_2$ and an unique function $\mathbf{\xi} : J_0 \times (-\epsilon_1, \epsilon_1) \to \mathbb{R}^2$ such that $\xi(\mu,0)=\mathbf{x} _\mu $ and $P(\xi(\mu, \epsilon), \mu, \epsilon) = \xi(\mu, \epsilon)$. This guarantees the existence of a periodic solution for \eqref{eqavrnormalform2}.

Consider the following hypothesis:

\begin{itemize}
\item[\bf{H2.}] Assume that $a' (\mu_0) \neq 0$.
\end{itemize}

In these conditions we have the following:

\begin{lemma}\label{lema} (\cite{candido2020torus}) Assume that hypotheses {\bf H1} and {\bf H2} hold. Let $\lambda(\mu, \epsilon)$ and $\overline{\lambda(\mu, \epsilon)}$ be the conjugate eigenvalues of $D_\mathbf{x} P(\xi(\mu, \epsilon), \mu, \epsilon)$.  Then there exists an unique smooth function $\mu: (-\epsilon_2,\epsilon_2) \to J_0$, with $0<\epsilon_2<\epsilon_1$ and $\mu(0)=\mu_0$ such that
\begin{itemize}
\item[(i)] $| \lambda(\mu(\epsilon), \epsilon) |= 1,$ 
\item[(ii)] $(\lambda(\mu(\epsilon), \epsilon) )^k \neq 1$, for $k=1,2,3,4$,
\item[(iii)] $\left.\dfrac{d}{d\mu} |\lambda(\mu,\epsilon)| \right|_{\mu=\mu(\epsilon)} \neq 0.$
\end{itemize}
\end{lemma}

Consider the change of variables $\mathbf{x} \rightarrow \mathbf{y} + \xi(\mu, \epsilon)$ and take $\mu = \sigma + \mu(\epsilon)$ to write the Poincar\'e map \eqref{poincaremap} as
\begin{equation}\label{poincaremap2}
\mathbf{y} \mapsto \mathbf{H}_\epsilon (\mathbf{y},\sigma) := \mathbf{y} + \epsilon \mathbf{G}(\mathbf{y} + \xi(\sigma + \mu(\epsilon), \epsilon), \sigma + \mu(\epsilon), \epsilon).
\end{equation}

With this expression, one can check that $\mathbf{y}=0$ is a fixed point for \eqref{poincaremap2}. Let $\eta_\epsilon(\sigma) = \lambda(\sigma + \mu(\epsilon), \epsilon)$ and $\overline{\eta_\epsilon(\sigma)} = \overline{\lambda(\sigma + \mu(\epsilon), \epsilon)}$ be the eigenvalues of $D_\mathbf{y} \mathbf{H}_\epsilon (0, \sigma)$. In its polar form, $\eta_\epsilon(\sigma) = r_\epsilon(\sigma) \text{e}^{i \phi_\epsilon(\sigma)}$, where $\phi_\epsilon(0) = \theta_\epsilon$, with $0 < \theta_\epsilon < \pi$. Properties \emph{(i)-(iii)} above translate into
$$
r_\epsilon(0)=1, \quad \text{e}^{i k \theta_\epsilon}\neq 1 \text{ for }k=1,2,3,4, \quad \text{ and } r'_\epsilon (0) \neq 0.
$$

The Taylor expansion of $\mathbf{H}_\epsilon(\mathbf{y},0)$ around $\mathbf{y}=0$ is
$$
\mathbf{H}_\epsilon(\mathbf{y},0) = \mathbf{A}_\epsilon \mathbf{y} + \frac{1}{2} \mathbf{B}_\epsilon (\mathbf{y},\mathbf{y}) + \frac{1}{6} \mathbf{C}_\epsilon (\mathbf{y},\mathbf{y},\mathbf{y})+\mathcal{O}(|\mathbf{y}|^4),
$$
where  $\mathbf{A}_\epsilon = D_\mathbf{y}\mathbf{H}_\epsilon(0,0)= Id +\epsilon \mathbf{\mathcal{A}}_\epsilon + \mathcal{O}(\epsilon^3)$ and $\mathbf{B}(\mathbf{u},\mathbf{v})$, $\mathbf{C}(\mathbf{u},\mathbf{v}, \mathbf{w})$ are multilinear functions that can be written as  $\mathbf{B}_\epsilon = \epsilon \mathcal{B}_\epsilon+ \mathcal{O}(\epsilon^3)$ and $\mathbf{C}_\epsilon = \epsilon \mathcal{C}_\epsilon + \mathcal{O}(\epsilon^3)$, where
$$
\begin{array}{rcl}
\mathcal{A}_\epsilon &=& A_1 + \epsilon A_2,\\ \\
\mathcal{B}_\epsilon (\mathbf{u},\mathbf{v}) &=& B_1(\mathbf{u},\mathbf{v}) + \epsilon B_2(\mathbf{u},\mathbf{v}),\\ \\
\mathcal{C}_\epsilon(\mathbf{u},\mathbf{v}, \mathbf{w}) &=& C_1(\mathbf{u},\mathbf{v}, \mathbf{w}) + \epsilon C_2(\mathbf{u},\mathbf{v}, \mathbf{w}).
\end{array}
$$

We can assume that $Id + \epsilon \mathcal{A}_\epsilon$ is in its real Jordan form without any loss of generality. This implies that
$$
Id + \epsilon \mathcal{A}_\epsilon = \begin{pmatrix} 1 + \alpha_\epsilon & -\beta_\epsilon \\ \beta_\epsilon & 1 + \alpha _\epsilon \end{pmatrix},
$$
where $\alpha_\epsilon \pm i \beta_\epsilon$ are the eigenvalues of $\mathcal{A}_\epsilon$.

This allows us to compute the first Lyapunov coefficient of the map \eqref{poincaremap2} at $(\mathbf{x},\sigma)=(0,0)$. Define $\mathbf{p}= \left( \frac{1}{\sqrt{2}}, \frac{-i}{\sqrt{2}} \right)$. Consider the usual inner product in $\mathbb{C}^2$ and let

\begin{equation}\label{lyapunov}
\begin{array}{rcl}
\ell_1^\epsilon &=& -\text{Re} \left( \dfrac{(1 - 2 \text{e}^{i \theta_\epsilon)})\text{e}^{-2 i \theta_\epsilon}}{2 - 2 \text{e}^{i\theta_\epsilon}} \left\langle \mathbf{p}, \epsilon \mathcal{B}_\epsilon (\mathbf{p},\mathbf{p}) \right\rangle \left\langle \mathbf{p}, \epsilon \mathcal{B}_\epsilon (\mathbf{p},\overline{\mathbf{p}}) \right\rangle \right)  \\ \\
&& + \text{Re}\left( \text{e}^{- i \theta_\epsilon} \left\langle \mathbf{p}, \epsilon \mathcal{C}_\epsilon (\mathbf{p},\mathbf{p},\overline{\mathbf{p}}) \right\rangle\right) - \dfrac{ \vert \left\langle \mathbf{p}, \epsilon \mathcal{B}_\epsilon (\mathbf{p},\overline{\mathbf{p}}) \right\rangle\vert ^2}{2}- \dfrac{ \vert \left\langle \mathbf{p}, \epsilon \mathcal{B}_\epsilon (\overline{\mathbf{p}},\overline{\mathbf{p}}) \right\rangle\vert ^2}{4}.
\end{array}
\end{equation}
This number can be written as a power series in $\epsilon$ as
$$
\ell_1^{\epsilon} = \epsilon \ell_{1,1} + \epsilon^2 \ell_{1,2} + \mathcal{O}(\epsilon^3).
$$

We state below the result used to prove Theorem \ref{teoremaprincipal} (b).

\begin{theorem}\label{teoremamurilo} (\cite{candido2020torus})
Assume that hypotheses {\bf H1} and {\bf H2} hold and let $j \in \{1,2\}$ be the smallest index such that and that $\ell_{1,j} \neq 0$. Then, for each $\epsilon>0$ sufficiently small there exists a smooth curve $\mu(\epsilon) \in J_0$ with $\mu(0)=\mu_0$, neighborhoods $U_\epsilon \subset \mathbb{S}^1 \times \Omega$ of the periodic solution $\phi(t, \mu(\epsilon),\epsilon)$ and intervals $J_\epsilon \subset J_0$ of $\mu(\epsilon)$ for which the following are true:
\begin{itemize}
\item[(i)] For $\mu \in J_\epsilon$ such that $\ell_1 (\mu - \mu(\epsilon)) \geq0$, the periodic solution $\phi(t,\mu(\epsilon), \epsilon)$ is unstable (resp. asymptotically stable) provided that $\ell_1 >0$ (resp. $\ell_1<0)$ and the differential equation \eqref{eqavrnormalform2} does not admit any invariant tori in $U_\epsilon$.
\item[(ii)] For $\mu \in J_\epsilon$ such that $\ell_1 (\mu - \mu(\epsilon)) <0$, the differential equation \eqref{eqavrnormalform2} admits a unique invariant torus $T_{\mu,\epsilon}$ in $U_\epsilon$ surrounding the periodic solution $\phi(t, \mu, \epsilon)$. Moreover, $T_{\mu,\epsilon}$ is unstable (resp. asymptotically stable) and the periodic solution $\phi(t,\mu,\epsilon)$ is asymptotically stable (resp. unstable).
\item[(iii)] $T_{\mu,\epsilon}$ is the unique torus of the differential equation \eqref{eqavrnormalform2} bifurcating from the periodic solution $\phi(t, \mu(\epsilon),\epsilon)$ in $U_\epsilon$ as $\mu$ passes through $\mu(\epsilon)$.
\end{itemize}

\end{theorem}

In Section 3 we verify hypotheses {\bf H1}-{\bf H2} and compute the coefficient $\ell_1$ explicitly for system \eqref{langfordsystem}.

\section{Proof of Theorems  \ref{teoremaprincipal0} and \ref{teoremaprincipal}}

\subsection{Proof of Theorem \ref{teoremaprincipal0}}

The usual averaging method is applied to prove item (a). Recall that under the perturbation \eqref{perturbation}, system \eqref{langfordsystem} has a zero-Hopf equilibrium at the origin. We put the linear form of this system in its Jordan normal form under the change of variables, along with a re-escaling by
$$
(x,y,z) = ( \epsilon Y,-\epsilon X,\epsilon Z).
$$

In these coordinates our system becomes
$$
\begin{array}{rcl}
\dot{X} &=& -Y \omega +\epsilon  \left(-\alpha _1 X+\mu _1 X+X Z-\beta _1 Y\right) +\epsilon ^2 \left(-\alpha _2 X+\mu _2 X-\beta _2 Y\right) + \mathcal{O}(\epsilon^3), \\ \\
\dot{Y} &=& X \omega +\epsilon  \left(\beta _1 X-\alpha _1 Y+Y \left(\mu _1+Z\right)\right) +\epsilon ^2 \left(\beta _2 X-\alpha _2 Y+\mu _2 Y\right) \mathcal{O}(\epsilon^3),\\ \\
\dot{Z} &=& \epsilon  \left(\mu_1 Z -\gamma _0 \left(X^2+Y^2+Z^2\right)\right) +\epsilon ^2 \left(\mu _2 Z -\gamma _1 \left(X^2+Y^2+Z^2\right)\right) + \mathcal{O}(\epsilon^3).
\end{array}
$$

Switching to cylindrical coordinates: $(X,Y,Z) = (r \cos (\theta), r \sin (\theta),z)$, one notices that $\dot{\theta} = \omega + \mathcal{O}(\epsilon)$. Since $\omega$ and $\epsilon$ are positive, then $\dot{\theta}>0$ and, for $\epsilon$ sufficiently small, one can use $\theta$ as the new time, obtaining

\begin{equation}\label{sistemareduzido}
\begin{array}{rcl}
\dfrac{dr}{d\theta} &=& \epsilon\dfrac{r \left(-\alpha _1+\mu _1+z\right)}{\omega } + \epsilon ^2 r  \left(\dfrac{\mu _2-\alpha _2}{\omega }-\dfrac{\beta _1 \left(-\alpha _1+\mu _1+z\right)}{\omega ^2}\right) + \mathcal{O}(\epsilon^3),\\ \\
\dfrac{dz}{d\theta} &=& -\epsilon \dfrac{  \left(\gamma _0 \left(r^2+z^2\right)-\mu _1 z\right)}{\omega } \\ \\
&&  +\epsilon ^2 \left(\dfrac{\beta _1 \left(\gamma _0 \left(r^2+z^2\right)-\mu _1 z\right)}{\omega ^2}-\dfrac{\gamma _1 \left(r^2+z^2\right)-\mu _2 z}{\omega }\right) + \mathcal{O}(\epsilon^3).
\end{array}
\end{equation}

We remark that this system is on the normal form \eqref{eqavrnormalform} and is autonomous in $\theta$ up to order 2 in $\epsilon$. Therefore, the averaging procedure is trivial at first and second order and the averaged functions are obtained directly. Using $T=2\pi$, from \eqref{g1geral} and \eqref{g2geral} we obtain
\begin{equation}\label{g1}
\mathbf{g}_1(r,z) = \big(\mathbf{g}_1^1(r,z),\mathbf{g}_1^2(r,z)\big)=\left( \dfrac{r \left(-\alpha _1+\mu _1+z\right)}{\omega }, -\dfrac{  \left(\gamma _0 \left(r^2+z^2\right)-\mu _1 z\right)}{\omega } \right)
\end{equation}
and
\begin{equation}\label{g2}
\mathbf{g}_2(r,z) = \big(\mathbf{g}_2^1(r,z),\mathbf{g}_2^2(r,z)\big)
\end{equation}
where
$$
\begin{array}{rcl}
\mathbf{g}_2^1(r,z) &=&- \dfrac{2 \pi  r}{\omega ^2} \left( -\pi  \alpha _1^2+\alpha _1 \left(2 \pi  \left(\mu _1+z\right)-\beta _1\right)+\beta _1 \left(\mu _1+z\right)
\right. \\ \\
 && \left. \omega  \left(\alpha _2-\mu _2\right)+\pi  \left(\gamma _0 \left(r^2+z^2\right)-z^2-\mu _1 \left(\mu _1+3 z\right)\right)
 \right)
\end{array}
$$ 
and
$$
\begin{array}{rcl}
\mathbf{g}_2^2(r,z) &=& \dfrac{2\pi}{\omega^2} \left( 2 \pi  \gamma _0^2 z \left(r^2+z^2\right)-\gamma _1 \omega  \left(r^2+z^2\right)+\gamma _0\left(-2 \pi  r^2 z\right. 
\right.\\ \\
&& \left.\left.+2 \pi  \alpha _1 r^2+\left(\beta _1-3 \pi  \mu _1\right) \left(r^2+z^2\right)\right)+z \left(-\beta _1 \mu _1+\mu _2 \omega +\pi  \mu _1^2\right) \right)
\end{array}
$$
are the first-order and second-order averaged functions of \eqref{sistemareduzido}, respectively. 

The non-linear system $\mathbf{g}_1(r,z)=(0,0)$ has only one solution on the domain of cylindrical coordinates $r>0$, which is
\begin{equation}\label{peqg1}
(r^*,z^*)=\left( \sqrt{\frac{\left(\mu _1-\alpha _1\right) \left(\alpha _1 \gamma _0-\left(\gamma _0+1\right) \mu _1\right)}{\gamma _0}}, \alpha _1-\mu _1 \right).
\end{equation}

From conditions \eqref{conditions} this solution is well defined. The Jacobian determinant of $\mathbf{g}_1$ at $(r^*,z^*)$ is
$$
\det D \mathbf{g}_1 (r^*,z^*) = \frac{2 \left(\mu_1 -\alpha _1\right) \left(\alpha _1 \gamma _0-\left(\gamma _0+1\right) \mu _1\right)}{\omega ^2},
$$
which, from conditions \eqref{conditions}, is non-zero. 

Moreover, the eigenvalues of $D \mathbf{g}_1 (r^*,z^*)$ are
\begin{equation}\label{eigenvalues}
\begin{array}{rcl}
\lambda_\pm&=&\dfrac{2 \gamma _0 \left(\mu _1-\alpha _1\right)+\mu _1}{2 \omega} \\ \\
&&  \pm \dfrac{\sqrt{4 \alpha _1^2 \gamma _0 \left(\gamma _0+2\right)-4 \alpha _1 \left(\gamma _0+2\right) \left(2 \gamma _0+1\right) \mu _1+\left(2 \gamma _0+3\right){}^2 \mu _1^2}}{2 \omega}.
\end{array}
\end{equation}
From \eqref{conditions}, these eigenvalues are complex conjugates. Proposition \ref{zerodapromed} guarantees the existence of a unique isolated periodic solution  $\phi(t,\epsilon)$ satisfying $\phi\rightarrow 0$ as $\epsilon \rightarrow 0$. The stability of this periodic solution is determined by the sign of $2 \gamma _0 \left(\mu _1-\alpha _1\right)+\mu _1$. This finishes the proof of item (a).

\subsection{Proof of item (b)}

Without any loss of generality, we assume $\alpha_1>0$. The case $\alpha_1<0$ is similar. Recall the first-order averaged system \eqref{g1}, taking into account the parameter $\mu_1$ and using $\mathbf{y} = (r,z)$:
\begin{equation}\label{sistemaaveragedmu}
\dot{\mathbf{y}}=\epsilon \mathbf{g}_1(\mathbf{y};\mu_1).
\end{equation}

First, we need to check Hypotheses $\mathbf{H1}$ and $\mathbf{H2}$. The eigenvalues of $D\mathbf{g}_1(\mathbf{y};\mu_1)$ are $a(\mu_1) \pm \,i \, b(\mu_1)$ given by \eqref{eigenvalues}, which have real part
\begin{equation}\label{partereal}
a(\mu_1) = \dfrac{2 \gamma _0 (\mu_1 - \alpha_1)+\mu _1}{2 \omega }.
\end{equation}
Since over the bifurcation parameter $\mu_1 = \widehat{\mu_0}$ this real part is zero, we get
\begin{equation}\label{mu0b}
\widehat{\mu_0} = \frac{2 \alpha _1 \gamma _0}{2 \gamma _0+1}
\end{equation}
and also from the imaginary part $b(\mu_1)$ in \eqref{eigenvalues},
\begin{equation}\label{omega0b}
\omega_0 = b(\widehat{\mu_0}) = \dfrac{\alpha _1 \sqrt{2\gamma _0}}{2 \gamma _0 \omega +\omega } >0.
\end{equation}
Note that, if $\alpha_1 <0$ (as mentioned above), one can take $-\omega_0$ instead of $\omega_0$ as the defining expression above and the proof follows similarly from here on.

	When $\mu_1=\widehat{\mu_0}$, $(\mathbf{y}_{\widehat{\mu_0}};\widehat{\mu_0})=(r^*,z^*)$ is a zero-Hopf equilibrium of system \eqref{sistemaaveragedmu}. This verifies Hypothesis $\mathbf{H1}$. From \eqref{partereal} one easily checks that
$$
a'(\widehat{\mu_0})=\dfrac{2 \gamma _0+1}{2 \omega }
$$
which is non-zero from conditions \eqref{conditions}, and thus $\mathbf{H2}$ is also verified.

Let us compute $\ell_1$. The first and second order averaged functions were obtained in the previous section, and the Poincar\'e map is
	\begin{equation}
	\mathbf{x}\mapsto P(\mathbf{x},\epsilon,\mu)  = \mathbf{x} + \epsilon \mathbf{G}( \mathbf{x}, \mu, \epsilon)
	\end{equation}
where
\begin{equation}
\mathbf{G}(\mathbf{x},\mu,\epsilon) = \mathbf{g}_1(\mathbf{x},\mu) + \epsilon \mathbf{g}_2 (\mathbf{x},\mu) + \epsilon^2 \widetilde{\mathbf{G}}(\mathbf{x},\mu,\epsilon).
\end{equation}

 For this, we use $\mu=\mu_1$ as the secondary bifurcation parameter. In \cite{candido2020torus} there is a closed-form expression for $\ell_{1,j}$, but to avoid lenghty calculations using the general case we work with power series using the averaged functions $\mathbf{g}_1$ and $\mathbf{g}_2$ previously obtained. We start using the Taylor series for $\xi(\mu, \epsilon)$ in $\epsilon$ around zero up to order 2. Assume that
\begin{equation}\label{ximueps}
\xi(\mu, \epsilon) = (a_0,b_0) + \epsilon (a_1,b_1) + \epsilon^2 (a_2,b_2) + \mathcal{O}(\epsilon^3).
\end{equation}
Since $\xi(\mu, \epsilon)$ is a fixed point of \eqref{poincaremap}, then it is a zero for \eqref{poincaremapG}. Substituting \eqref{ximueps} in \eqref{poincaremapG}, one can solve $\mathbf{G}(\xi(\mu_1, \epsilon), \mu_1, \epsilon)=(0,0)$ for $a_0,b_0,a_1,b_1$, obtaining:
$$
\begin{array}{rcl}
a_0 &=& \sqrt{\dfrac{\left(\mu_1-\alpha _1\right) \left(\alpha _1 \gamma _0-\left(\gamma _0+1\right) \mu _1\right)}{\gamma _0}},\\ \\
b_0 &=&  \alpha _1-\mu _1,  \\ \\
a_1 &=& \dfrac{\pi}{\gamma _0^{3/2} \sqrt{\left(\alpha_1-\mu_1\right) \left(\gamma _0 \left(\mu_1 -\alpha _1\right)+\mu_1 \right)}} \left( \alpha _2 \gamma _0 \left(2 \gamma _0 \left(\mu_1 -\alpha _1\right)+\mu_1 \right)\right. \\ \\
&& \left.+\gamma _1 \mu_1  \left(\mu_1 -\alpha _1\right)+\gamma _0 \mu _2 \left(2 \gamma _0 \left(\alpha _1-\mu_1 \right)+\alpha _1-2 \mu_1 \right) \right),  \\ \\
b_1 &=&  2 \pi  \left(\alpha _2-\mu _2\right), \\ \\
a_2 &=&   -\dfrac{\pi^2}{2 \gamma _0^{5/2} \omega  \left(\left(\alpha _1-\mu_1 \right) \left(\gamma _0 \left(\mu_1 -\alpha _1\right)+\mu_1 \right)\right){}^{3/2}} \left(2 \alpha _2 \gamma _0 \left(\gamma _0 \mu_1  \left(2 \mu_1  \left(\mu_1 -\alpha _1\right) \left(-2 \pi  \alpha _1\right.\right.\right.\right.\\ \\
&& \left.\left.\left.\left. +\beta _1+2 \pi  \mu_1 \right)-\alpha _1 \mu _2 \omega \right)+4 \gamma _0^3 \left(\mu_1 -\alpha _1\right){}^3 \left(-\pi  \alpha _1+\beta _1+\pi  \mu_1 \right)\right.\right.\\ \\
&& \left.\left. +2 \gamma _0^2 \mu_1  \left(\mu_1 -\alpha _1\right){}^2 \left(-4 \pi  \alpha _1+3 \beta _1+4 \pi  \mu_1 \right)+\gamma _1 \mu_1 ^2 \omega  \left(\alpha _1-\mu_1 \right)\right)\right.\\ \\
&& \left. +\gamma _1 \mu_1  \left(\mu_1 -\alpha _1\right){}^2 \left(4 \beta _1 \gamma _0 \left(\gamma _0 \left(\mu_1 -\alpha _1\right)+\mu_1 \right)+\gamma _1 \omega  \left(4 \gamma _0 \left(\alpha _1-\mu_1 \right)-3 \mu_1 \right)\right)\right.\\ \\
&& \left. -2 \gamma _0 \mu _2 \left(\mu_1 -\alpha _1\right) \left(2 \gamma _0 \left(\gamma _0 \left(\mu_1 -\alpha _1\right)+\mu_1 \right) \left(-\alpha _1 \left(2 \gamma _0+1\right) \left(\beta _1+2 \pi  \mu_1 \right)\right.\right.\right.\\ \\
&& \left. \left.\left. +2 \pi  \alpha _1^2 \gamma _0+2 \left(\gamma _0+1\right) \mu_1  \left(\beta _1+\pi  \mu_1 \right)\right)+\gamma _1 \omega  \left(\mu_1  \left(\alpha _1-2 \mu_1 \right)-2 \gamma _0 \left(\mu_1 -\alpha _1\right){}^2\right)\right)\right.\\ \\
&& \left. +\alpha _2^2 \gamma _0^2 \mu_1 ^2 \omega +\alpha _1^2 \gamma _0^2 \mu _2^2 \omega \right),\\ \\
b_2 &=& \dfrac{4 \pi ^2}{\gamma _0 \omega} \left( \gamma _0 \left(\alpha _2 \beta _1-\mu _2 \left(-\pi  \alpha _1+\beta _1+\pi  \mu_1 \right)\right)+\pi  \gamma _1 \mu_1  \left(\mu_1 -\alpha _1\right)\right).
\end{array}
$$

Let $\lambda(\mu(\epsilon), \epsilon)$ and $\overline{\lambda(\mu(\epsilon), \epsilon)}$ be the eigenvalues of $D_x \mathbf{P} (\xi(\mu(\epsilon), \epsilon), \mu(\epsilon), \epsilon)$. Assume that their real and imaginary parts are
$$
\text{Re} \lambda(\mu(\epsilon), \epsilon) = 1+ \epsilon a + \epsilon ^2 b + \mathcal{O}(\epsilon^3),
$$
$$
\text{Im} \lambda(\mu(\epsilon), \epsilon) =\epsilon \omega_0 + \epsilon^2 \omega_1 + \mathcal{O}(\epsilon^3),
$$
for some $a,b, \omega_1 \, \in\mathbb{R}$, where $\omega_0$ is defined in \eqref{omega0}. Let $$\mu(\epsilon) = \widehat{\mu_0} +  \widehat{\mu_1}\epsilon +\widehat{\mu_2} \epsilon^2 + \mathcal{O}(\epsilon^3)$$ be the Taylor expansion of the curve from Lemma \ref{lema} around $0$ in $\epsilon$, where $\widehat{\mu_0}$ is obtained above. Then $\lambda(\mu(\epsilon), \epsilon)$ satisfies
\begin{itemize}
\item[\emph{(i)}] $| \lambda(\mu(\epsilon), \epsilon) |= 1,$ 
\item[\emph{(ii)}] $(\lambda(\mu(\epsilon), \epsilon) )^k \neq 1$, for $k=1,2,3,4$,
\item[\emph{(iii)}] $\left.\dfrac{d}{d\mu} |\lambda(\mu(\epsilon),\epsilon)| \right|_{\mu=\mu(\epsilon)} \neq 0.$
\end{itemize}

We remark that Lemma \ref{lema} is an existence result and it doesn't provide an explicit expression for $\mu(\epsilon)$. We can work with power series in $\epsilon$ to obtain it, which will satisfy items \emph{(i)} - \emph{(iii)} up to order 2 in $\epsilon$. In fact, we use these items as conditions to determine $a$ and $b$ in terms or $\omega_0$ and $\omega_1$ and then use the actual expressions of $\lambda(\mu(\epsilon), \epsilon)$ and its conjugate calculated using $P (\xi(\mu_1, \epsilon), \mu_1, \epsilon)$  to compute $\widehat{\mu_1}$ and $\widehat{\mu_2}$.

Item \emph{(i)} above can be used to deduce that
$$
a=0, \quad b= - \dfrac{\omega_0^2}{2}.
$$
With this, we have
$$
\begin{array}{rclll}
(\lambda(\mu(\epsilon), \epsilon) )^2 &=& 1+2 i \omega_0 \epsilon +\epsilon ^2 \left(-2 \omega_0^2+2 i \omega_1\right)+\mathcal{O}\left(\epsilon ^3\right) & \neq & 1,\\ \\
(\lambda(\mu(\epsilon), \epsilon) )^3 &=& 1+3 i \omega_0 \epsilon +\epsilon ^2 \left(-\frac{9 \omega_0^2}{2}+3 i \omega_1\right)+\mathcal{O}\left(\epsilon ^3\right) &\neq & 1, \\ \\
(\lambda(\mu(\epsilon), \epsilon) )^4 &=& 1+4 i \omega_0 \epsilon +\epsilon ^2 \left(-8 \omega_0^2+4 i \omega_1\right)+\mathcal{O}\left(\epsilon ^3\right)&\neq & 1.
\end{array}
$$

Now, explicitely, the real part of $\lambda(\mu(\epsilon), \epsilon)$ as a series of $\epsilon$ up to order 2 is
$$
\begin{array}{rcl}
\text{Re } \lambda(\mu(\epsilon), \epsilon) &=&1 + \epsilon^2\dfrac{1}{2 \left(2 \gamma _0 \omega +\omega \right){}^2}\left(8 \gamma _0^3 \omega  \left(-2 \pi  \alpha _2+2 \pi  \mu _2+\widehat{\mu}_1\right)\right.\\ \\
&& \left.+4 \gamma _0^2 \omega  \left(-4 \pi  \alpha _2+6 \pi  \mu _2+3 \widehat{\mu}_1\right)-2 \gamma _0 \left(4 \pi  \alpha _1 \gamma _1 \omega +2 \pi  \alpha _2 \omega \right. \right. \\ \\
&& \left. \left. +4 \pi ^2 \alpha _1^2-3 \omega  \left(2 \pi  \mu _2+\widehat{\mu}_1\right)\right)+\omega  \left(-4 \pi  \alpha _1 \gamma _1+2 \pi  \mu _2+\widehat{\mu}_1\right)\right) + \mathcal{O}(\epsilon^3).
\end{array}
$$
Since this equals to $1-\epsilon^2 \dfrac{\omega_0^2}{2}+ \mathcal{O}(\epsilon^3)$, we conclude that
$$
\begin{array}{rcl}
\widehat{\mu_1} &=& \dfrac{2}{\left(2 \gamma _0+1\right){}^3 \omega } \left(\pi  \left(2 \gamma _0+1\right){}^2 \omega  \left(2 \alpha _2 \gamma _0-\left(2 \gamma _0+1\right) \mu _2\right) \right.\\ \\
&& \left. +2 \pi  \alpha _1 \left(2 \gamma _0+1\right) \gamma _1 \omega +\left(4 \pi ^2-1\right) \alpha _1^2 \gamma _0\right).
\end{array}
$$

Using the imaginary part of $\lambda(\mu(\epsilon), \epsilon)$ as a series of $\epsilon$ up to order 2 and the expressions of $\widehat{\mu_0}$ and $\widehat{\mu_1}$, we obtain

$$
\begin{array}{rcl}
\omega_1 &=& \dfrac{\sqrt{2}}{\sqrt{\gamma _0} (2 \gamma _0+1){}^3 \omega ^2} \left(2 \pi  \gamma _0 \left(2 \gamma _0+1\right){}^2 \left(\alpha _2 \omega -\alpha _1 \beta _1\right) \right.\\ \\
&& \left. +\pi  \alpha _1 \left(1-4 \gamma _0^2\right) \gamma _1 \omega +4 \pi ^2 \alpha _1^2 \gamma _0-\alpha _1^2 \gamma _0\right).
\end{array}
$$

This allows us to compute
$$
\begin{array}{rcl}
\widehat{\mu_2} &=& \dfrac{4}{\left(2 \gamma _0+1\right){}^5 \omega ^2} \left( \pi ^2 \left(2 \gamma _0+1\right){}^3 \omega  \left(2 \alpha _2 \left(\beta _1 \gamma _0 \left(2 \gamma _0+1\right)+\gamma _1 \omega \right)-\beta _1 \left(2 \gamma _0+1\right){}^2 \mu _2\right) \right.\\ \\
&& \left. +\pi  \alpha _1^2 \left(2 \gamma _0+1\right) \left(\left(1+4 \pi ^2\right) \beta _1 \gamma _0 \left(2 \gamma _0+1\right)+\left(4 \left(1-6 \pi ^2\right) \gamma _0-1\right) \gamma _1 \omega \right) \right.\\\\ 
&& \left. +2 \pi  \alpha _1 \left(2 \gamma _0+1\right){}^2 \omega  \left(\left(4 \pi ^2-1\right) \alpha _2 \gamma _0+\pi  \gamma _1 \left(2 \beta _1 \gamma _0+\beta _1-2 \gamma _1 \omega \right)\right) \right.\\ \\
&&\left. +\left(1-4 \pi ^2\right)^2 \alpha _1^3 \gamma _0 \right).
\end{array}
$$
and

$$
\begin{array}{rcl}
\omega_2 &=& \dfrac{1}{2 \sqrt{2} \gamma_0^{3/2} \left(2 \gamma_0+1\right){}^5 \omega^3} \left(4 \pi \alpha_1^2 \gamma_0 \left(2 \gamma_0+1\right)\left(4 \beta_1 \gamma_0 \left(2 \gamma_0+1\right)\right.\right. \\ \\
&& \left.
+ \left(2 \gamma_0\left(16 \pi^2 \gamma_0 \left(\gamma_0+1\right)-24 \pi^2+5\right)-4 \pi^2-1\right)\gamma_1 \omega \right)\\ \\
&& +2 \pi \alpha_1 \left(2 \gamma_0+1\right){}^2 \omega\left(8 \alpha_2
\left(4 \pi^2 \gamma_0 \left(\gamma_0+1\right)+5 \pi^2-1\right)\gamma_0^2\right.\\ \\
&&\left.+\pi\left(\left(12 \left(\gamma_0-1\right)\gamma_0-1\right)\gamma_1^2 \omega-16 \pi \gamma_0^2\left(2 \gamma_0+1\right){}^2 \mu_2\right)\right)\\ \\
&&-8 \pi^2 \alpha_2 \gamma_0\left(2 \gamma_0-1\right)\left(2\gamma_0+1\right){}^3\gamma_1 \omega^2\\ \\
&& \left.+\left(4 \pi^2-1\right)\alpha_1^3\gamma_0^2\left(\gamma_0\left(4\gamma_0\left(\left(1+12\pi^2\right)\gamma_0+4 \pi^2+3\right)-20 \pi^2+9\right)+24 \pi^2-6\right)\right).
\end{array}
$$

Finally we have that
$$
\begin{array}{rcl}
\left.\dfrac{d}{d\mu_1} \vert \lambda (\mu, \epsilon) \vert  \right|_{\mu_1 = \mu(\epsilon)} &=& \dfrac{-1}{\sqrt{2\gamma _0} \omega } + \epsilon \dfrac{\alpha _1 \gamma _0 \left(-2 \pi ^2 \gamma _0+4 \pi ^2-1\right)+\pi  \left(2 \beta _1 \gamma _0+\gamma _1 \omega \right)}{\sqrt{2} \gamma _0^{3/2} \omega ^2} \\ \\
&&+ \mathcal{O}(\epsilon^2) \neq 0.
\end{array}
$$
This shows that, up to order $\epsilon^2$, the curve $\mu(\epsilon)$ explicitated here is the curve that exists and is unique due to Lemma \ref{lema}.

By performing the change of variables $\mathbf{x} \rightarrow \mathbf{y} + \xi(\mu, \epsilon)$ and taking $\mu = \sigma + \mu(\epsilon)$, we write the Poincar\'e map as
\begin{equation}\label{poincare3}
\mathbf{y} \mapsto \mathbf{H}_\epsilon (\mathbf{y},\sigma) = \mathbf{y} + \epsilon \mathbf{G}(\mathbf{y} + \xi(\sigma + \mu(\epsilon), \epsilon), \sigma + \mu(\epsilon), \epsilon).
\end{equation}
Note that $\mathbf{y}=0$ is a fixed point of \eqref{poincare3}. From the eigenvalues $\eta_\epsilon(\sigma)=\lambda(\sigma + \mu(\epsilon), \epsilon)$ and $\overline{\eta_\epsilon(\sigma)}=\overline{\lambda(\sigma + \mu(\epsilon), \epsilon)}$ of $D_\mathbf{y} \mathbf{H}_\epsilon (0, \sigma)$, we write

$$
\eta_\epsilon (\sigma) = r_\epsilon (\sigma) \, \text{e}^{i \phi_\epsilon (\sigma)},
$$
with $\phi_\epsilon (0) = \theta_\epsilon \in (0,\pi)$. From Lemma \ref{lema}, $r_\epsilon(0) = 1$ and $\text{e}^{i k \theta_\epsilon} \neq 1$ for $k=1,2,3,4$. In particular, we have
$$
\theta_\epsilon = \dfrac{\pi }{2}+\epsilon\dfrac{\omega_0}{2}-\epsilon ^2\dfrac{\omega_1}{2}+\mathcal{O}(\epsilon^3).
$$

Let the Taylor expansion of $\mathbf{H}_\epsilon(\mathbf{y},0)$ around $\mathbf{y}=0$ be
\begin{equation}\label{taylorpoincare}
\mathbf{H}_\epsilon(\mathbf{y},0) = \mathbf{A}_\epsilon \mathbf{y} + \frac{1}{2} \mathbf{B}_\epsilon (\mathbf{y},\mathbf{y}) + \frac{1}{6} \mathbf{C}_\epsilon (\mathbf{y},\mathbf{y},\mathbf{y})+\mathcal{O}(|\mathbf{y}|^4),
\end{equation}
where
$$
\mathbf{A}_\epsilon = \text{Id} + \epsilon \mathcal{A}_ \epsilon + \mathcal{O}(\epsilon^3)
$$
and
$$
\mathcal{A}_ \epsilon = A_1 + \epsilon A_2.
$$

We find that
$$
A_1 = \left(
\begin{array}{cc}
 0 & -\frac{\alpha _1}{2 \gamma _0 \omega +\omega } \\
 \frac{2 \alpha _1 \gamma _0}{2 \gamma _0 \omega +\omega } & 0 \\
\end{array}
\right)
$$
and
$$
A_2=\left(
\begin{array}{cc}
 -\frac{4 \pi ^2 \alpha _1^2 \gamma _0}{\left(2 \gamma _0 \omega +\omega \right){}^2} & \frac{m_1}{\left(2 \gamma _0+1\right){}^3 \omega ^2} \\
 \frac{m_2}{\left(2 \gamma _0+1\right){}^3 \omega ^2} & \frac{2 \left(2 \pi ^2-1\right) \alpha _1^2 \gamma _0}{\left(2 \gamma _0 \omega +\omega \right){}^2} \\
\end{array}
\right)
$$
where
$$\begin{array}{rcl}
m_1&=&2 \pi  \alpha _1 \left(2 \gamma _0+1\right) \left(\beta _1 \left(2 \gamma _0+1\right)+2 \gamma _1 \omega \right)-2 \pi  \alpha _2 \left(2 \gamma _0+1\right){}^2 \omega +\left(1-4 \pi ^2\right) \alpha _1^2
\end{array}
$$
and
$$
\begin{array}{rcl}
m_2&=&-4 \pi  \alpha _1 \left(2 \gamma _0+1\right) \left(\beta _1 \gamma _0 \left(2 \gamma _0+1\right)-\gamma _1 \omega \right) +4 \pi  \alpha _2 \gamma _0 \left(2 \gamma _0+1\right){}^2 \omega \\
&&+2 \left(4 \pi ^2-1\right) \alpha _1^2 \gamma _0 .
\end{array}
$$

Consider the change of variables $\mathbf{y} = T_\epsilon \mathbf{u}$, where $\mathbf{u} = (u,v)^T$,  $$T_\epsilon =  T_0 + \epsilon T_1 + \epsilon^2 T_2 + \mathcal{O}(\epsilon^3)$$ and
$$
T_0=\left(
\begin{array}{cc}
 \frac{1}{\sqrt{2\gamma _0}} & 0 \\
 0 & 1 \\
\end{array}
\right), \quad
 T_1 = \left(
\begin{array}{cc}
 \frac{1}{\sqrt{2\gamma _0} } & \frac{1}{\sqrt{2\gamma _0}} \\
 \frac{\left(1-4 \pi ^2\right) \alpha _1 \sqrt{\gamma _0}}{\sqrt{2} \left(2 \gamma _0 \omega +\omega \right)}-1 & \frac{\pi  \gamma _1}{\gamma _0}+1 \\
\end{array}
\right),
$$
and $T_2 = \begin{pmatrix} \frac{1}{\sqrt{2 \gamma _0} } & \frac{1}{\sqrt{2 \gamma _0} } \\ t_{1,1} & t_{1,2} \end{pmatrix},$ where
$$
\begin{array}{rcl}
t_{1,1} &=& \dfrac{\alpha _1 \left(2 \pi  \left(8 \pi ^2 \gamma _0+1\right) \gamma _1 \omega -\gamma _0 \left(2 \gamma _0+1\right) \left(2 \left(\pi +4 \pi ^3\right) \beta _1+\left(4 \pi ^2-1\right) \omega \right)\right)}{\sqrt{2} \sqrt{\gamma _0} \left(2 \gamma _0 \omega +\omega \right){}^2}\\
&& -\dfrac{\left(1-4 \pi ^2\right)^2 \alpha _1^2 \sqrt{\gamma _0}}{\sqrt{2} \left(2 \gamma _0+1\right){}^3 \omega ^2}+\dfrac{\pi  \left(1-4 \pi ^2\right) \sqrt{2} \alpha _2 \sqrt{\gamma _0}}{2 \gamma _0 \omega +\omega }-\dfrac{\pi  \gamma _1}{\gamma _0}-1,
\\
\\
t_{1,2} &=& \dfrac{1}{8} \left(-\dfrac{2 \left(1-4 \pi ^2\right)^2 \alpha _1^2 \gamma _0}{\left(2 \gamma _0 \omega +\omega \right){}^2} +\dfrac{4 \left(1-4 \pi ^2\right) \sqrt{2} \alpha _1 \sqrt{\gamma _0}}{2 \gamma _0 \omega +\omega } \right.\\ \\
&& \left. +\dfrac{4 \pi  \gamma _1 \left(2 \gamma _0 \left(2 \pi  \beta _1+\omega \right)-\pi  \gamma _1 \omega \right)}{\gamma _0^2 \omega }+8\right).\\
\end{array}
$$

On the new coordinates $(u,v)$, $\mathcal{A}_\epsilon$ is its real Jordan form which is
$$
\mathcal{A}_\epsilon = \begin{pmatrix}
1 - \epsilon^2 \dfrac{\omega_0^2}{2} & -\epsilon \omega_0 + \epsilon^2 \omega_1 \\
\epsilon \omega_0 - \epsilon^2 \omega_1& 1 - \epsilon^2 \dfrac{\omega_0^2}{2}
\end{pmatrix} + \mathcal{O}(\epsilon^3).
$$

Let $\mathbf{u}=(u_1,u_2)$, $\mathbf{v} = (v_1,v_2)$ and $\mathbf{w} = (w_1,w_2)$. With $\mathcal{A}_\epsilon$ as above and putting \eqref{taylorpoincare} on coordinates $(u,v)$, we obtain  $\mathcal{B}_\epsilon = B_0 + \epsilon B_1$ and $\mathcal{C_\epsilon} = C_0 + \epsilon C_1$.  We have:
$$
B_0=\left(\dfrac{u_1 v_2+u_2 v_1}{2 \omega },-\dfrac{u_1 u_2+2 \gamma _0 v_1 v_2}{2 \omega } \right),
$$
$B_1 = (B_1^1, B_1^2)$, where
$$
\begin{array}{rcl}
B_1^1 &=& \dfrac{1}{2 \gamma _0 \omega ^2} \left(\dfrac{\sqrt{2} \alpha _1 \gamma _0^{3/2} \left(u_1 \left(2 \pi ^2 u_2+u_2\right)+4 \pi ^2 \gamma _0 v_1 v_2-4 \pi ^2 v_1 v_2\right)}{2 \gamma _0+1} \right. \\ \\
&& +\gamma _0 \left(-2 \pi  \beta _1 (u_1 v_2+u_2 v_1)   +u_1 \omega  (v_2-u_2)+v_1 \omega  (u_2+2 v_2)\right)\\ \\
&&\left. +\pi  \gamma _1 \omega  (u_1 v_2+u_2 v_1)+2 \gamma _0^2 v_1 v_2 \omega \right),
\\ \\
B_1^2 &=& \dfrac{1}{4 \omega ^2}\left(\dfrac{\sqrt{2} \alpha _1 \sqrt{\gamma _0} \left(-2 \gamma _0+12 \pi ^2-1\right) (u_1 v_2+u_2 v_1)}{2 \gamma _0+1}+4 \pi  \beta _1 \left(u_1 u_2+2 \gamma _0 v_1 v_2\right) \right.\\ \\
&& \left.+4 \gamma _0 \omega  (v_2 (u_1-v_1)+u_2 v_1)-\dfrac{2 \pi  \gamma _1 \omega  \left(u_1 u_2+6 \gamma _0 v_1 v_2\right)}{\gamma _0}-2 u_1 u_2 \omega \right).
\end{array}
$$

We also obtain $C_0 = (0,0)$, $C_1 = (C_1^1,C_1^2)$, which are
$$
\begin{array}{rcl}
C_1^1 &=& \dfrac{\pi ^2}{3 \omega ^2} \left(-2 \gamma _0 (u_1 v_2 w_2+u_2 v_1 w_2+w_1 v_1 v_2)-3 u_1 u_2 w_1+2 u_1 v_2 w_2+2 u_2 v_1 w_2+2 w_1 v_1 v_2\right), \\ \\
C_1^2 &=&\dfrac{2 \pi ^2}{3 \omega ^2} \left(\gamma _0 \left(u_1 u_2 w_2+u_1 w_1 v_2+u_2 w_1 v_1+6 \gamma _0 v_1 v_2 w_2\right)-u_2 (u_1 w_2+w_1 v_1)-u_1 w_1 v_2\right).
\end{array}
$$

Finally, let $\mathbf{p} = \left( \frac{1}{\sqrt{2}}, -\frac{i}{\sqrt{2}}\right) \in \mathbb{C}^2$. We use \eqref{lyapunov} to obtain
$$
\ell_1^\epsilon = \dfrac{\left(4 \pi ^2-2\right) \gamma _0^2+3 \left(4 \pi ^2-1\right) \gamma _0+\pi ^2}{16 \omega ^2} \epsilon^2 + \mathcal{O}(\epsilon^3).
$$
Under the assumption that $\ell_{1,2} \neq 0$ the proof is finished. As a final remark, since $\gamma_0>0$ it is straightforward to see that $\ell_{1,2} >0$, and, from theorem \ref{teoremamurilo}, no asymptotically stable torus can bifurcate from the periodic solution.

\section{Numerical Examples}

In this section we provide numerical examples for the existence of periodic orbits and invariant tori in system \eqref{langfordsystem} using the conditions established in Theorem \ref{teoremaprincipal}.

\begin{example}\label{exemplo1} \textit{An isolated asymptotically stable periodic solution.}

\begin{figure}[h!] 
\centering
\includegraphics[scale=0.4]{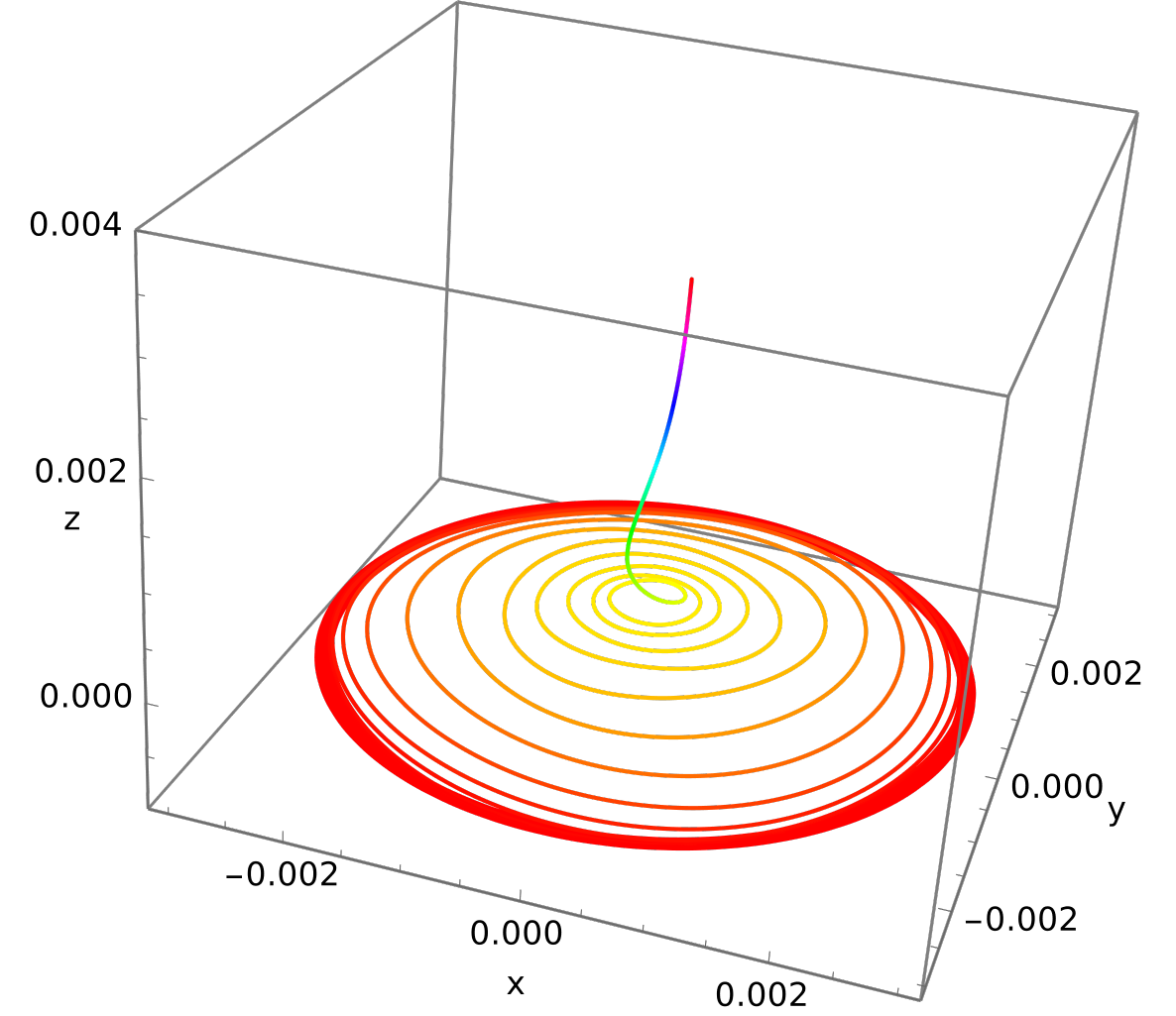}
\caption{The trajectory starting at $(0.0001, 0.0001, 0.0010)$ being attracted to the periodic solution from Example \ref{exemplo1}.}
\label{fig1}
\end{figure}

Consider system \eqref{langfordsystem} with parameters
\begin{align*}
\alpha = - \epsilon, \quad
\beta = \dfrac{1}{100} , \\
\gamma = \dfrac{2}{5}- \dfrac{1}{20}\epsilon, \quad
\mu  =-\dfrac{11}{12} \epsilon.
\end{align*}

From Theorem \ref{teoremaprincipal0}, with these values system \eqref{langfordsystem} has a periodic solution near the origin; this periodic solution is asymptotically stable since
$$
2 \gamma _0 \left(\mu _1-\alpha _1\right)+\mu _1 = -\dfrac{17}{20} <0.
$$
We remark that no invariant torus bifurcates from the periodic solution. Figure \ref{fig1} illustrates the behaviour of a solution near the origin for these values.
\end{example}

\begin{example}\label{exemplo2} \textit{An invariant unstable torus bifurcating from a periodic solution.}

Consider system \eqref{langfordsystem} with parameters
\begin{align*}
\alpha = \epsilon + \dfrac{1}{100} \epsilon^2, \quad
\beta = \dfrac{1}{10}+ \dfrac{443}{100}\epsilon + \dfrac{1}{25} \epsilon^2 , \\
\gamma = -1- \epsilon+   \dfrac{1}{25} \epsilon^2, \quad
\mu  = \dfrac{399317}{100000} \epsilon+  \dfrac{1}{25} \epsilon^2.
\end{align*}

With these values we have
$$
\ell_{1,2} = \dfrac{49 \pi ^2-11}{144} >0 .
$$

Using a small $\epsilon$ such as $\epsilon = \frac{1}{1000}$ we have $\mu_1 - \mu(\epsilon) \approx -7.4 \cdot 10^{-3}$. Theorem \ref{teoremamurilo} guarantees the existence of a periodic orbit of \eqref{langfordsystem} surrounded by an unstable invariant torus. In Figure \ref{fig2} we see the orbit with initial conditions $(0.0005, 0, 0.0005)$. In Figure \ref{fig3} the Poincar\'e map of this orbit is illustrated.

\begin{figure}[H] 
\centering
\includegraphics[scale=0.6]{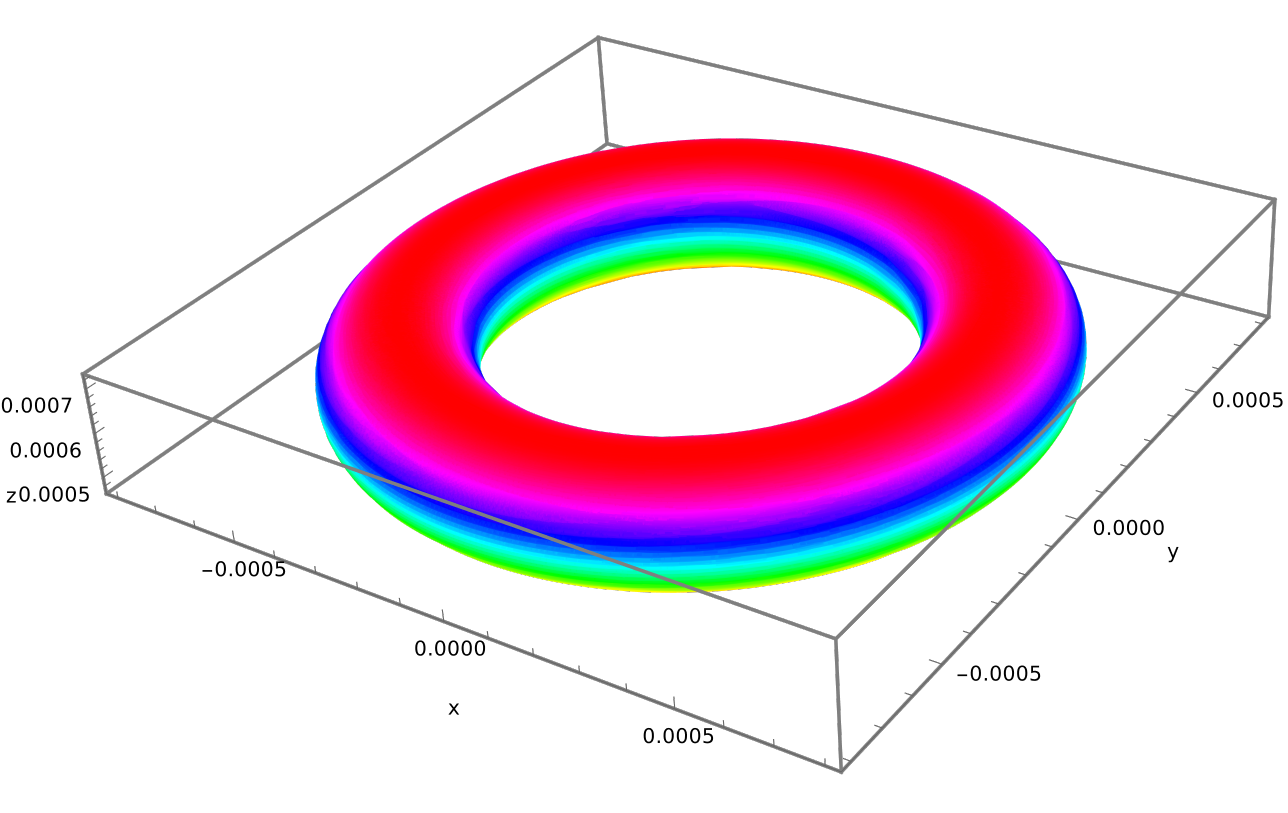}
\caption{The trajectory starting at $(0.0005, 0, 0.0005)$ from Example \ref{exemplo2}.}
\label{fig2}
\end{figure}

\begin{figure}[H] 
\centering
\includegraphics[scale=0.22]{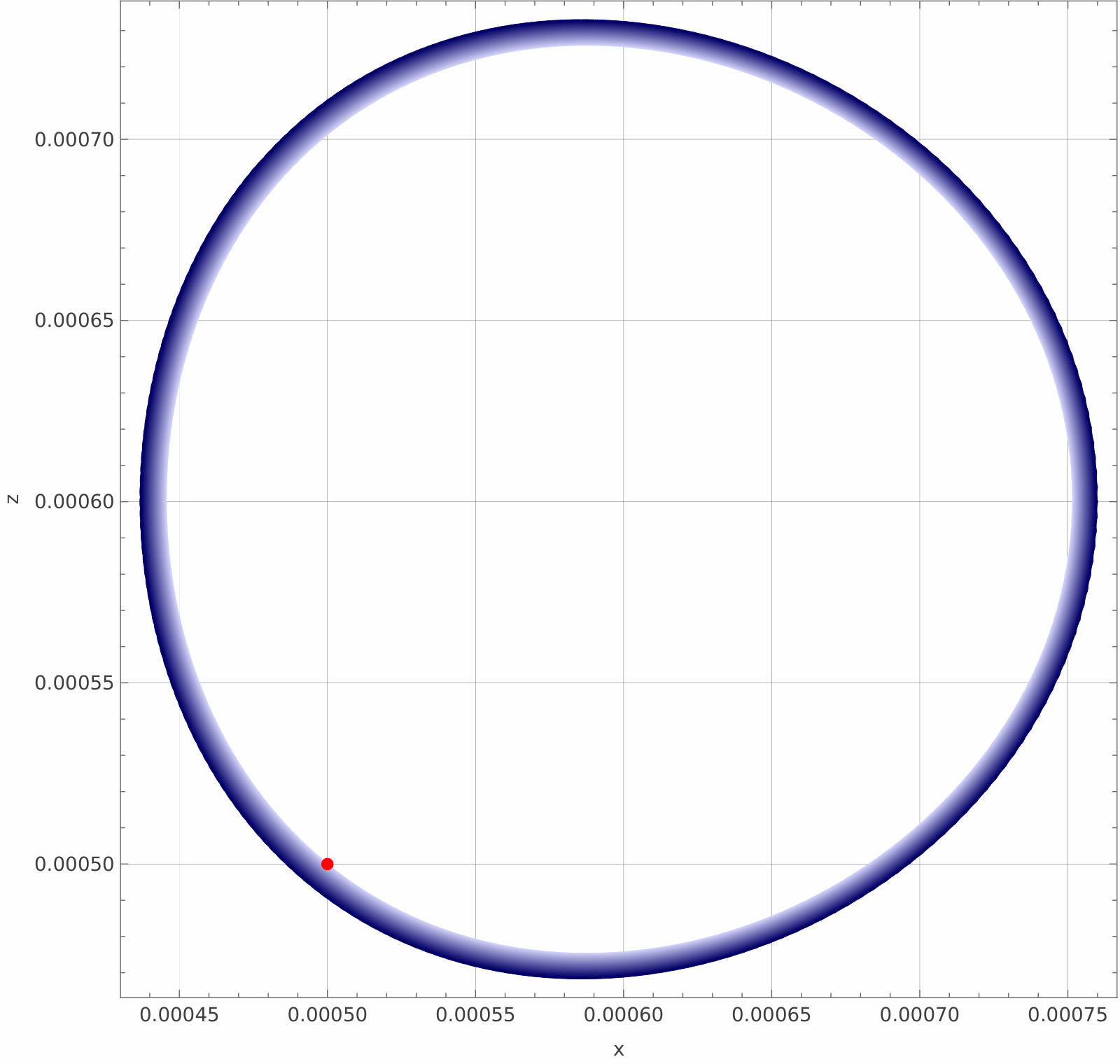}
\caption{Poincar\'e map of the trajectory considered in Example~\ref{exemplo2}, 
starting at $(x_0, y_0, z_0) = (0.0005, 0, 0.0005)$, obtained from the 
intersections of the orbit with the half-plane $y = 0$, $x > 0$. The 
red dot marks the first such intersection (which coincides with the 
initial condition, since $y_0 = 0$); subsequent intersections are 
shown in progressively darker shades of blue as time increases.}
\label{fig3}
\end{figure}

\end{example}

\bmhead{Acknowledgements}

I would like to express my gratitude to Murilo R. C\^andido for introducing me to this fascinating field of research and  for his invaluable help throughout the development of this work.

\section*{Declarations}

\begin{itemize}
\item \textbf{Funding:} The author did not receive support from any organization for the submitted work.
\item \textbf{Competing Interests:} The author has no competing interests to declare that are relevant to the content of this article.
\item \textbf{Ethics approval and consent to participate:} Not applicable.
\item \textbf{Consent for publication:} Not applicable.
\item \textbf{Data availability:} No datasets were generated or analysed during the current study; all computations reported can be reproduced from the equations presented in the manuscript.
\item \textbf{Materials availability:} Not applicable.
\item \textbf{Code availability:} The symbolic computations and numerical simulations were performed using Wolfram Mathematica. The corresponding notebooks are available from the author upon reasonable request.
\item \textbf{Author contribution:} The author confirms sole responsibility for the following: study conception, analytical derivations, numerical simulations, and manuscript preparation.
\end{itemize}


\end{document}